\documentclass[12pt]{amsart}
\usepackage{amssymb}
\def\NZQ{\mathbb}               

\def\CC{{\NZQ C}}

\def\opn#1#2{\def#1{\operatorname{#2}}} 
\opn\chara{char}
\opn\length{\ell}
\opn\projdim{proj\,dim}
\opn\injdim{inj\,dim}
\opn\rank{rank}
\opn\depth{depth}
\opn\grade{grade}
\opn\height{height}
\opn\embdim{emb\,dim}
\opn\codepth{codepth}
\opn\codim{codim}

\opn\Tr{Tr}
\opn\bigrank{big\,rank}
\opn\superheight{superheight}\opn\lcm{lcm}
\opn\trdeg{tr\,deg}%
\opn\reg{reg}
\opn\ini{in}
\opn\div{div}
\opn\Div{Div}
\opn\cl{cl}
\opn\Cl{Cl}
\opn\Spec{Spec}
\opn\Supp{Supp}
\opn\supp{supp}
\opn\Sing{Sing}
\opn\Ass{Ass}
\opn\Ann{Ann}
\opn\Rad{Rad}
\opn\Soc{Soc}
\opn\Ker{Ker}
\opn\Coker{Coker}
\opn\Im{Im}
\opn\Hom{Hom}
\opn\Tor{Tor}
\opn\Ext{Ext}
\opn\End{End}
\opn\Aut{Aut}
\opn\id{id}

\opn\nat{nat}
\opn\pff{pf}
\opn\Pf{Pf}
\opn\GL{GL}
\opn\SL{SL}
\opn\mod{mod}
\opn\ord{ord}
\opn\aff{aff}
\opn\con{conv}
\opn\relint{relint}
\opn\st{st}
\opn\lk{lk}
\opn\cn{cn}
\opn\core{core}
\opn\vol{vol}
\opn\gr{gr}

\def\pot#1#2{#1[\kern-0.28ex[#2]\kern-0.28ex]}

\opn\dirlim{\underrightarrow{\lim}}
\opn\invlim{\underleftarrow{\lim}}
\let\to=\rightarrow

\def\Implies{\ifmmode\Longrightarrow \else
     \unskip${}\Longrightarrow{}$\ignorespaces\fi}
\def\implies{\ifmmode\Rightarrow \else
     \unskip${}\Rightarrow{}$\ignorespaces\fi}
\def\iff{\ifmmode\Longleftrightarrow \else
     \unskip${}\Longleftrightarrow{}$\ignorespaces\fi}

\let\:=\colon
\newtheorem{Theorem}{Theorem}[section]
\newtheorem{Lemma}[Theorem]{Lemma}
\newtheorem{Corollary}[Theorem]{Corollary}

\newtheorem{Definition}[Theorem]{Definition}

\let\epsilon=\varepsilon
\let\phi=\varphi
\let\kappa=\varkappa
\def\qed{\ifhmode\textqed\fi
   \ifmmode\ifinner\quad\qedsymbol\else\dispqed\fi\fi}
\def\textqed{\unskip\nobreak\penalty50
    \hskip2em\hbox{}\nobreak\hfil\qedsymbol
    \parfillskip=0pt \finalhyphendemerits=0}
\def\dispqed{\rlap{\qquad\qedsymbol}}

\begin{document}
\title[Bivariate Composed Product]{A Bivariate Analogue to the Composed\\ Product of Polynomials}
\author{Donald Mills \& Kent M. Neuerburg}
\address{Donald Mills\\ Department of Mathematics, Southern Illinois University, Carbondale, IL 62901}
\email{dmills@math.siu.edu}
\address{Kent M. Neuerburg, Mathematics Department, Southeastern Louisiana University, Hammond, LA 70402}
\email{kneuerburg@selu.edu}
\thanks {The first author was a Davies Fellow for the National Research Council.  He wishes to thank the NRC, and specifically the U.S. Army Research Laboratory and the U.S. Military Academy, for the use of their facilities during the time that this paper was completed.}
\begin{abstract}
The concept of a {\it composed product} for univariate polynomials has been explored extensively by Brawley, Brown, Carlitz, Gao, Mills, et al.  Starting with these fundamental ideas and utilizing fractional power series representation (in particular, the Puiseux expansion) of bivariate polynomials, we generalize the univariate results.  We define a bivariate composed sum, composed multiplication, and composed product (based on function composition).  Further, we investigate the algebraic structure of certain classes of bivariate polynomials under these operations.  We also generalize a result of Brawley and Carlitz concerning the decomposition of polynomials into irreducibles.
\end{abstract}
\maketitle

\section{The Univariate Composed Product}
To begin, let $\mathbb{F}_{q}$ denote the finite field of $q$ elements, $q$ a prime power, let $\mathbb{F}_{q}[x]$ denote the integral domain of polynomials in the indeterminate $x$ over $\mathbb{F}_{q}$, and let $\Gamma_{q}$ denote the algebraic closure of $\mathbb{F}_{q}$. Further, let $G$ denote a nonempty subset of $\Gamma_{q}$ which enjoys the following properties:

\begin{enumerate}
\item{$G$ is invariant under the Frobenius automorphism $\sigma$ where $\sigma\hspace{0.02in}:\hspace{0.02in}\alpha\mapsto\alpha^{q}$ for $\alpha \in \Gamma_{q}$; and}
\item{There is defined on $G$ a binary operation $\diamond$ such that for all $\alpha$, $\beta \in G$, $\sigma(\alpha \diamond \beta)=\sigma(\alpha) \diamond \sigma(\beta)$.}
\end{enumerate}
These properties will be referred to as the \emph {$\diamond$-properties}.

In the notation of Brawley and Carlitz \cite{bc1}, we denote by $M_{G}[q,x]$ the set of all nonconstant monic polynomials $f \in \mathbb{F}_{q}[x]$ whose roots lie in $G$. Then the {\it composed product\/} of the polynomials $f$, $g \in \mathbb{F}_{q}[x]$ is defined as

\begin{eqnarray}
\displaystyle (f \odot g)(x)=\prod_{\alpha}\prod_{\beta}(x-\alpha \diamond \beta),
\end{eqnarray}
where the ordinary products $\prod$ are over all roots $\alpha$ of $f$ and $\beta$ of $g$, including multiplicites, so that the degree of $f \odot g$ is the product of the degrees of $f$ and $g$. (Remark: The notation $\diamond$ has been used in much of the literature to denote both the operation on the group $G$ as well as the binary operation on $M_{G}[q,x]$. In an effort to avoid confusion as well as provide a more natural framework for what is to follow, we have decided to use the $\odot$ operation to denote the composed product operation on $M_{G}[q,x]$.)

The composed product operation (or $\odot$-operation) is a binary operation on $M_{G}[q,x]$ \cite{bc1}, as can be seen by observing that $[(f \odot g)(x)]^q=(f \odot g)(x^q)$.  What is also true is that the $\odot$-operation inherits many of the properties of the binary operation $\diamond$ which induces the $\odot$-operation. For instance, if $G$ is a semigroup under $\diamond$ then $M_{G}[q,x]$ is a semigroup under $\odot$; if $G$ is commutative under $\diamond$ then $M_{G}[q,x]$ is commutative under the $\odot$-operation. If $G$ possesses an identity element $e \in G \cap \mathbb{F}_{q}$, then the polynomial $x-e$ serves as the identity relative to the $\odot$-operation. Finally, it is clear that the units of $M_{G}[q,x]$ are the polynomials $x-u$ where $u \in G \cap \mathbb{F}_{q}$ is a unit in $G$.

In their seminal work, Brawley and Carlitz deal with an issue which is analogous to ordinary factorization of polynomials, namely the {\it decomposition} of a polynomial in $M_{G}[q,x]$ according to the composed product, where $G$ is now an abelian group under $\diamond$. We define the following. A polynomial $f \in M_{G}[q,x]$ which is not a unit is said to be {\it decomposable} with respect to the $\odot$-operation if and only if there exist polynomials $g$, $h \in M_{G}[q,x]$, each of degree at least two, such that $f(x)=(g \odot h)(x)$. Otherwise, $f$ is said to be {\it indecomposable\/}. Polynomials $f$, $g \in M_{G}[q,x]$ are said to be {\it associates} if and only if $f(x)=(r \odot g)(x)$ for some unit $r \in M_{G}[q,x]$, and in this case we write $f \sim g$. Clearly, $\sim$ is an equivalence relation on $M_{G}[q,x]$.

Before launching into a detailed discussion of decomposition of polynomials according to $\odot$, Brawley and Carlitz prove the following theorem which serves to narrow one's focus with regards to the decomposition question.
\begin{Theorem}
Suppose that $(G,\diamond)$ is a group and let $f$, $g \in M_{G}[q,x]$ be given with deg$(f)=m$ and deg$(g)=n$. Then the composed product $f \odot g$ is irreducible over $\mathbb{F}_{q}$ if and only if $f$ and $g$ are both irreducible over $\mathbb{F}_{q}$ with $\gcd(m,n)=1$.
\end{Theorem}

Theorem 1.1 encourages us to ask whether irreducibles in $M_{G}[q,x]$ decompose uniquely, up to associates, into irreducible indecomposables. Brawley and Carlitz show in \cite{bc1} that the answer to this question is yes when the binary operation on the group $G$ is either ordinary multiplication or ordinary addition, and along the way they provide tests to determine whether a given polynomial is decomposable with respect to either the composed multiplication operation \cite{bc1}, the term used when the operation on $G$ is ordinary multiplication, or the composed sum operation \cite{bc2}, the term used when the operation on $G$ is ordinary addition. Brawley and Brown generalized the work in \cite{bc1} to include all abelian groups $G$ \cite{bb}. The general statement and its corollary are given below.

\begin{Theorem}\label{T:decomp1}
Let $G \subset \Gamma_{q}$ be a group with binary operation $\diamond$ which satisfies the $\diamond$-properties listed above, let $f(x) \in M_{G}[q,x]$ be irreducible over $\mathbb{F}_{q}$ with deg$(f)=n>1$, and let $\odot$ denote the composed product operation on $M_{G}[q,x]$. Suppose that $f$ can be decomposed in two ways as 
\begin{center}
$f(x)=(f_{1} \odot f_{2} \odot \cdots \odot f_{t})(x)=(g_{1} \odot g_{2} \odot \cdots \odot g_{s})(x)$,
\end{center}
where each of the $f_{i}$ and $g_{i}$ belong to $M_{G}[q,x]$ and are each indecomposable with respect to the $\odot$-operation. Then $s=t$ and there is some reordering of the $g_{i}$'s so that $f_{i} \sim g_{i}$ for $i=1,2,...,t$.
\end{Theorem}
We say that $a \in G$ is {\it indecomposable} with respect to $\diamond$ ($G$'s binary operation) if the minimal polynomial of $a$ over $\mathbb{F}_q$ cannot be decomposed nontrivially with respect to the composed product operation $\odot$.  This definition provides the following corollary to Theorem \ref{T:decomp1}.
\begin{Corollary} Let $G \subset \Gamma_{q}$ be a group with binary operation $\diamond$ which satisfies the properties listed above. Every element $\gamma \in G \setminus \mathbb{F}_{q}$ can be written as a $\diamond$-product $\gamma=\alpha_{1} \diamond \alpha_{2} \diamond \cdots \diamond \alpha_{t}$ of a finite number of indecomposables, where deg$(\gamma)$ over $\mathbb{F}_{q}$ equals $\prod_{i=1}^{t}$deg$(\alpha_{i})$. Moreover, if $\gamma=\beta_{1} \diamond \beta_{2} \diamond \cdots \diamond \beta_{s}$ is another such $\diamond$-decomposition of $\gamma$, then $s=t$ and there is a reordering of the $\beta_{i}$'s such that for each $i$ from 1 to $t$, there exists a $c_{i} \in G \cap \mathbb{F}_{q}$ such that $\alpha_{i}=c_{i} \diamond \beta_{i}$, and further $c_{1} \diamond c_{2} \diamond \cdots \diamond c_{t}$ equals the identity element of $G$.
\end{Corollary}

Other issues that have been considered in the study of the composed product include the matter of simultaneous decomposition \cite{bc1}, in which it is asked whether a given irreducible can decompose nontrivially (i.e., none of the components in the decomposition is a unit) according to both the composed sum and composed multiplication operations (the answer is no); the efficient computation of a general form of the composed product, and in particular the efficient computation of the composed multiplication and composed sum of two polynomials \cite{bgm1}; the factorization patterns of reducible polynomials $h \in M_{G}[q,x]$ which decompose as 
$h(x)=(f \odot g)(x)$ for irreducibles $f$, $g \in M_{G}[q,x]$ of non-coprime degrees \cite{mil}; and the determination of the group structure of $G$ when the binary operation on $G$ is represented by a bivariate rational function $R(x,y) \in k(x,y) \setminus k[x,y]$, where $k(x,y)$ is the function field in the variables $x$ and $y$ over a field $k$ and $k[x,y]$ is the domain of polynomials in $x$ and $y$ over $k$ \cite{bgm2}.

Our goal in this paper is to define a bivariate analogue to the $\odot$-operation described above. That is, given a field $k$, indeterminates $x$ and $y$, and polynomials $f(x,y)$, $g(x,y) \in k[x,y]$ which satisfy certain requirements, we use $f$ and $g$ to create a new polynomial $h(x,y) \in k[x,y]$ in a manner that is similar to the way in which we formed the composed product of two univariate polynomials, as given by (1). 

\section{Fractional Power Series.}  
Our approach to the bivariate case will be to avail ourselves of the well-known theorem of Puiseux (though an earlier version of the theorem was known to Newton).
\begin{Theorem}\label{T:puiseux}
Let $k$ be an algebraically closed field and let $f(x,y)=y^m+a_1(x)y^{m-1}+ \cdots +a_m(x) \in k((x))[y]$ be a monic polynomial with deg$_yf(x,y)=m>0$ and coefficients $a_1(x), \dots, a_m(x)$ in $k((x))$.  Further, assume either the characteristic of $k$ is zero or that $m!$ is not divisible by the characteristic of $k$.  Then there exists a positive integer $n$, not divisible by the characteristic of k, such that
$$f(x,y)=\prod\limits_{i=1}^m(y-p_i(x^{\frac{1}{n}}))$$ 
with $p_i(t) \in k((t))$.
\end{Theorem}
\noindent
Additionally, we have the following corollaries (see \cite{Ab}).
\begin{Corollary}
If $a_i(x) \in k[[x]]$ for $i=1,2, \dots, m$ then $p_i(t) \in k[[t]]$ for $i=1,2, \dots, m$.
\end{Corollary}
\begin{Corollary}\label{C:conjugates}
If $f(x,y)$ is irreducible in $k((x))[y]$ then we have $n=m$ (in fact $m$ is the least possible value of $n$).  Moreover, the power series $p_1(x), \dots p_n(x)$ may be arranged so that $p_i(x^{\frac{1}{n}})=p_1(\omega^i x^{\frac{1}{n}})$ where $\omega$ is a primitive $m^{th}$-root of unity.
\end{Corollary}

In letters to Oldenburg, Newton \cite{Ne} develops a constructive method for determining a $p_i(x)$.  This construction forms the basis of the proof of the general theorem, see \cite{Wa}.  It is important to note that Puiseux's theorem is valid in all characteristics.  For a more detailed discussion of the characteristic $p$ case see \cite{Ca}.

\section{The Bivariate Composed Product}
\subsection{Definitions}  
We are now prepared to generalize the univariate case. 
For $f(x,y), g(x,y) \in k((x))[y]$, we invoke Theorem \ref{T:puiseux} to write $f(x,y)=\prod\limits_{i=1}^{m_1}(y-p_i(x^{\frac{1}{n_1}}))$ and $g(x,y)=\prod\limits_{j=1}^{m_2}(y-q_j(x^{\frac{1}{n_2}}))$.  We now make the following definitions.
\begin{Definition}
The {\bf composed sum} of $f(x,y)$ and $g(x,y)$ is given by
$$
f \star g = \prod\limits_{i=1}^{m_1}\prod\limits_{j=1}^{m_2}(y-(p_i(x^{\frac{1}{n_1}})+q_j(x^{\frac{1}{n_2}})))
$$
\end{Definition}
\begin{Definition}
The {\bf composed multiplication} of $f(x,y)$ and $g(x,y)$ is given by
$$
f \bullet g = \prod\limits_{i=1}^{m_1}\prod\limits_{j=1}^{m_2}(y-(p_i(x^{\frac{1}{n_1}})q_j(x^{\frac{1}{n_2}})))
$$
\end{Definition}
\noindent
In the bivariate case we have an additional operation available to us, namely function composition. \begin{Definition}
If $f(0,0)=g(0,0)=0$ ($f$ and $g$ have $0$ as constant term) then let the {\bf composed product} of $f(x,y)$ and $g(x,y)$ be given by
$$
f \odot g = \prod\limits_{i=1}^{m_1}\prod\limits_{j=1}^{m_2}(y-p_i(q_j(x^{\frac{1}{n_2}})^{\frac{1}{n_1}}))
$$
\end{Definition}

We quickly illustrate each of these operations with an example (for computational convenience, we assume we are working over $\CC$).  Let $f(x,y)= y^4-2x^3y^2-4x^5y+x^6-x^7$.  Then, applying Newton's construction method we obtain $y=p(x^{\frac{1}{4}})=x^{\frac{6}{4}}+x^{\frac{7}{4}}$
as one solution (see \cite{BK}).  We observe that $f$ is irreducible so by Corollary \ref{C:conjugates}, we have $p_i(x^{\frac{1}{4}})=(\omega_1^ix^{\frac{1}{4}})^6+(\omega_1^ix^{\frac{1}{4}})^7=p(\omega_1^ix^{\frac{1}{4}})$ for $i=1, \dots, 4$, where $\omega_1$ is a primitive fourth root of unity.  Hence, 
$$
f(x,y)=\prod\limits_{i=1}^4(y-p(\omega_1^ix^{\frac{1}{4}})).
$$

Similarly, let $g(x,y)=y^6-3x^3y^4-2x^5y^3+3x^6y^2-6x^8y-x^9+x^{10}$.  Like $f$, $g$ is irreducible.  In this case, Newton's construction method yields a root $y=q(x^{\frac{1}{6}})=x^{\frac{9}{6}}+x^{\frac{10}{6}}$.  Hence, $q_j(x^{\frac{1}{6}})=(\omega_2^jx^{\frac{1}{6}})^9+(\omega_2^jx^{\frac{1}{6}})^{10}=q(\omega_2^jx^{\frac{1}{6}})$ for $j=1, \dots, 6$.  Here, $\omega_2$ is a primitive sixth root of unity.  Thus, 
$$
g(x,y)=\prod\limits_{j=1}^6(y-q(\omega_2^jx^{\frac{1}{6}})).
$$

Using these factorizations we compute the composed sum to be
\begin{equation*}
\begin{aligned}
f \star g&=\prod\limits_{i=1}^4\prod\limits_{j=1}^6(y-(p(\omega_1^ix^{\frac{1}{4}})+ q(\omega_2^jx^{\frac{1}{6}})))\\
	&=\prod\limits_{i=1}^4\prod\limits_{j=1}^6(y-(\omega_1^{6i}x^{\frac{6}{4}}+\omega_1^{7i}x^{\frac{7}{4}}+\omega_2^{9j}x^{\frac{9}{6}}+\omega_2^{10j}x^{\frac{10}{6}})).
\end{aligned}
\end{equation*}

A similar computation shows the composed multiplication is
\begin{equation*}
\begin{aligned}
f \bullet g&=\prod\limits_{i=1}^4\prod\limits_{j=1}^6(y-(p(\omega_1^ix^{\frac{1}{4}})q(\omega_2^jx^{\frac{1}{6}})))\\
	&=\prod\limits_{i=1}^4\prod\limits_{j=1}^6(y-(\omega_1^{6i}\omega_2^{9j}x^3 + \omega_1^{7i}\omega_2^{9j}x^{\frac{13}{4}} + \omega_1^{6i}\omega_2^{10j}x^{\frac{19}{6}} + \omega_1^{7i}\omega_2^{10j}x^{\frac{41}{12}}))
\end{aligned}
\end{equation*}

The composed product is, on the other hand, a much more complicated computation.  Then, using our example $f$ and $g$ we get (after expanding the terms as power series and simplifying)
\begin{equation*}
\begin{aligned}
f \odot g&=\prod\limits_{i=1}^4\prod\limits_{j=1}^6(y-p(\omega_1^iq(\omega_2^jx^{\frac{1}{6}})^{\frac{1}{4}}))\\
	&=\prod\limits_{i=1}^4\prod\limits_{j=1}^6(y-(\omega^{12i+6j}x^{\frac{54}{24}}+\frac{3}{2}\omega^{12i+10j}x^{\frac{58}{24}}+\frac{3}{8}\omega^{12i+14j}x^{\frac{62}{24}}+\omega^{18i+15j}x^{\frac{63}{24}}+\cdots))
\end{aligned}
\end{equation*}
where $\omega$ represents a primitive $24^{th}$ root of unity.

\subsection{Properties}
In what follows, $\mathbb{M}$ will represent the set of all polynomials in $k((x))[y]$ which satisfy the hypotheses of Theorem \ref{T:puiseux}, while $\mathbb{M}^{\ast} \subset \mathbb{M}$ will represent the set of all elements of $\mathbb{M}$ whose Puiseux expansions have roots $p_{i}$ which belong to $k((x)) \setminus k$ for all $i$.

It is natural to ask what sort of structure $\mathbb{M}$ possesses under the $\star$, $\bullet$, and $\odot$ operations. When the characteristic of $k$ is zero, we have the following.

\begin{Theorem}\label{T:structure1}
(a) $(\mathbb{M},\star)$ is a commutative semigroup with identity $e(x,y)=y$. (b) $(\mathbb{M},\bullet)$ is a commutative semigroup with identity $e(x,y)=y-1$. (c) $(\mathbb{M}^{\ast},\odot)$ is a non-commutative semigroup with identity $e(x,y)=y-x$.
\end{Theorem}
\begin{proof} We prove (a), as the proofs of (b) and (c) are similar. Select any pair of elements $f(x,y)$, $g(x,y) \in \mathbb{M}$, with deg$_{y}f(x,y)=m$ and deg$_{y}g(x,y)=n$; then clearly $(f \star g)(x,y) \in \mathbb{M}$ as $f \star g$ is monic in $y$ with deg$_{y}(f \star g)=mn>0$ and with coefficients in $k((x))$, and thus closure under $\star$ is proved. The associativity and commutativity of the $\star$ operation follow from the associativity and commutativity properties which $k((x))$ enjoys under ordinary addition. The identity element of $(\mathbb{M},\star)$ is clearly $e(x,y)=y$, while the identity elements for $(\mathbb{M},\bullet)$ and $(\mathbb{M},\odot)$ are $e(x,y)=y-1$ and $e(x,y)=y-x$ respectively.
\end{proof}

Note that in general we cannot make the same statements when the characteristic is $p>0$, for we can always find $f(x,y)$, $g(x,y) \in \mathbb{M}$ whose degrees in $y$ are $m$ and $n$ respectively such that $p \mid (mn)!$. However, if one considers the set of all elements of $\mathbb{M}$ whose degree in $y$ is $1$ (we denote this by $\mathbb{M}_{1}$), as well as the set of elements of $\mathbb{M}_{1}$ of the form $f(x,y)=y-a(x)$ with $a(x) \in k((x)) \setminus k$ (denoted by $\mathbb{M}_{1}^{\ast}$), we have the following.

\begin{Theorem}\label{T:structure2}
Let $\mathbb{M}_{1}$ be defined as above with char$(k) \ge 0$. Then $(\mathbb{M}_{1},\star,\bullet)$ is a field with $\star$-identity $e(x,y)=y$ and $\bullet$-identity $i(x,y)=y-1$.  Further, $(\mathbb{M}_1,\star,\bullet) \cong (k((x)),+,\cdot)$.
\end{Theorem}
\begin{proof} The proof that $(\mathbb{M}_1,\star,\bullet)$ is a field proceeds in the same fashion as the proof of Theorem \ref{T:structure1}.  

To see that $(\mathbb{M}_1,\star,\bullet) \cong (k((x)),+,\cdot)$ we consider the canonical map $\phi:\mathbb{M}_1 \to k((x))$ via $\phi: (y-a(x)) \mapsto a(x)$.  Note that if $f(x,y)=y-a(x)$ and $g(x,y)=y-b(x)$ then  $\phi(f(x,y) \star g(x,y))=\phi(y-(a(x)+b(x)))=a(x)+b(x)=\phi(f(x,y))+\phi(g(x,y))$.  Similarly, $\phi(f(x,y) \bullet g(x,y))=\phi(f(x,y))\cdot\phi(g(x,y))$.  Thus $\phi$ is a homomorphism of fields.  If $\phi(y-a(x))= \phi(y-b(x))$ then, clearly, $a(x)=b(x)$; hence, $\phi$ is injective.  For any $a(x) \in k((x))$ we have $a(x)=\phi(y-a(x))$, so $\phi$ is also surjective.
\end{proof}

For the $\odot$-operation, we have the following.
\begin{Theorem}
$(\mathbb{M}_1^{\ast},\odot)$ is a non-commutative semigroup with identity $e(x,y)=y-x$.
\end{Theorem}
\begin{proof}  Proceeds in the same fashion as the proof of Theorem \ref{T:structure1}.
\end{proof}

For certain subsets of $\mathbb{M}^{\ast}$ under $\odot$ with $char(k)=p>0$, a somewhat greater structure can be obtained. Specifically, let $\mathbb{M}_{h,n}$ denote the set of all elements of $\mathbb{M}^{\ast}$ which are homogeneous polynomials of degree $0<n<\sqrt{p}$ and whose $x^{n}$ and $y^n$ coefficients are nonzero (the reason for the second requirement will become obvious shortly), and set $\mathbb{M}_{h}=\bigcup_{n=1}^{\mathbb{[}\sqrt{p}\mathbb{]}}\mathbb{M}_{h,n}$, where $[s]$ is the greatest integer less than or equal to $s$. Using Newton's expansion method for determining the Puiseux series, it is easy to show that for any $f(x,y)=\sum_{i+j=n}a_{i,j}x^{i}y^{j} \in \mathbb{M}_{h}$, a Puiseux expansion of $f$ is given by $y=\alpha x$ where $\alpha$ is a root of the polynomial $w_{f}(t)=a_{0,n}t^{n}+a_{1,n-1}t^{n-1}+ \cdots +a_{n-1,1}t+a_{n,0}$, $a_{n,0} \neq 0$ (so the zero element of $k$ is not a root of $w_{f}$; this will be important when we employ the $\odot$-operation on a certain subset of $\mathbb{M}_{h}$). We call $w_{f}$ the {\it polynomial associated with $f$\/}. 

Let $\mathbb{M}_{h,min}$ denote the set of elements $f(x,y) \in \mathbb{M}_h$
whose associated polynomials $w_f$ act as their own minimum polynomials.
That is, if $a \in k$, $r(t) \in k[t]$ with $w_f(a)=0=r(a)$ then $w_{f}(t) \mid r(t)$ in $k[t]$. Thus, $w_{f}(t)$ is irreducible over some subfield of $k$; specifically, $w_{f}(t)$ is irreducible over $L$ where $L \subset k$ is the smallest subfield containing the coefficients of $w_{f}(t)$. Now let $k=\Gamma_{q}$ and let $f$, $g \in \mathbb{M}_{h,min}$ with $deg(f)=m$, $deg(g)=n$, $(m,n)=1$, and the coefficients of $f$ and $g$ belonging to $\mathbb{F}_{q}$.  Further let the associated polynomials of $f$ and $g$ be $w_f$ and 
$w_g$, respectively, with roots $a$ and $b$ respectively.  Then $F(x,y)=(f
\odot g)(x,y)$ is a homogeneous polynomial of degree $mn<p$ with associated
polynomial $w_F$ having root $ab$.  Since $(m,n)=1$ with $w_f$ and $w_g$ irreducible we
see by Theorem 1.1 that $w_F=w_f \odot w_g$ (where $(G,\diamond)$ is a
multiplicative group here) is an irreducible univariate polynomial of
degree $mn$ over $\mathbb{F}_q$ with root $ab$. Thus, $(f \odot g)(x,y) \in
\mathbb{M}_{h,min}$ and we have the following result.

\begin{Lemma}\label{L:well-defined} 
Let $k$ be an algebraically closed field having characteristic $p$, and let $f(x,y)$, $g(x,y) \in \mathbb{M}_{h,min}$ be homogeneous polynomials whose associated polynomials $w_{f}(t)$, $w_{g}(t) \in k[t]$ have coprime degrees. Then the $\odot$-operation, when performed on $f$ and $g$, is a well-defined operation, that is $F(x,y)=(f \odot g)(x,y)$ is a homogeneous polynomial of degree deg$(f)$deg$(g)$ whose associated polynomial $w_{F}$ has degree deg$(w_{f})$deg$(w_{g})$. Moreover, $(f \odot g) \odot h=f \odot (g \odot h)$ for all triplets $f(x,y)$, $g(x,y)$, $h(x,y) \in M_{h,min}$ whose associated polynomials have coprime degrees, and also $f \odot g=g \odot f$ for all pairs $f(x,y)$, $g(x,y) \in \mathbb{M}_{h,min}$ whose associated polynomials have coprime degrees. The polynomial $e(x,y)=y-x$ serves as the identity element relative to the $\odot$-operation.
\end{Lemma}
>From Lemma \ref{L:well-defined} we have the following.

\begin{Theorem}\label{T:decomp2}
Let $k$ be an algebraically closed field having characteristic $p$, and let $\mathbb{M}_{h,min,1} \subset \mathbb{M}_{h,min}$ denote the set of elements of $\mathbb{M}_{h,min}$ of degree 1. Then $(\mathbb{M}_{h,min,1},\odot)$ is an abelian group.  Further, $(\mathbb{M}_{h,min,1},\odot)\cong(k^*,\cdot)$.
\end{Theorem}
\begin{proof}
Select any pair of elements $f(x,y)$, $g(x,y) \in \mathbb{M}_{h,min,1}$, and write their Puiseux expansions as $y=ax$ and $y=bx$ respectively (where $a$, $b \neq 0$ and are each of degree 1 over $k$). Under the $\odot$-operation, the Puiseux expansion for the composed product of $f$ and $g$ is $y=(ab)x$.  Thus $(f \odot g)(x,y)$ has associated polynomial $w(t)$ which is of degree one, and hence is the minimum polynomial of $ab$ over $k$. The $\odot$-operation is clearly associative and commutative on $\mathbb{M}_{h,min,1}$, and it is clear as well that the identity element of $\mathbb{M}_{h,min,1}$ under $\odot$ is given by $e(x,y)=y-x$. As for inverses, we see that if $f(x,y)$, $g(x,y) \in \mathbb{M}_{h,min}$ have Puiseux expansions $y=ax$ and $y=a^{-1}x$ respectively then $(f \odot g)(x,y)=e(x,y)$. 

If $f(x,y) \in \mathbb{M}_{h,min,1}$ then $f(x,y)=y-ax$, $a\in k^*$.  Let $\phi: \mathbb{M}_{h,min,1} \to k^*$ by $\phi:(y-ax) \mapsto a$.  Let $f(x,y)=y-ax$ and $g(x,y)=y-bx$.  We have $\phi(f(x) \odot g(x))=\phi(y-abx)=ab=\phi(f(x))\cdot\phi(g(x))$ showing $\phi$ is a homomorphism.  If $\phi(f(x,y))=\phi(y-ax)=\phi(y-bx)=\phi(g(x,y))$ then $a=b$ so $f(x,y)=g(x,y)$ making $\phi$ injective.  For every $a\in k^*$ we have $a=\phi(y-ax)$, hence $\phi$ is surjective.  Thus $\phi$ is an isomorphism.
\end{proof}

What is particularly nice about the elements of $\mathbb{M}_{h,min}$ is that their definition, along with the definition of the $\odot$-operation applied to these polynomials, gives rise to a bivariate analogue of Theorem \ref{T:decomp1} in the case where the coefficients of the elements of $\mathbb{M}_{h,min}$ come from $\mathbb{F}_{q}((x))\subset\Gamma_q((x))$. To see this, we first adapt the definition of associates given in Section 2. Specifically, two polynomials $f(x,y)$, $g(x,y) \in \mathbb{M}_{h,min}$ are {\it associates\/} of each other under the $\odot$-operation provided there exists a polynomial $h(x,y) =y-ax \in \mathbb{F}_q[x,y]$ such that $f(x,y)=(h \odot g)(x,y)$. To denote this relationship, we write $f \sim g$, and note that this relation is an equivalence relation on $\mathbb{M}_{h,min}$. We have the following.

\begin{Theorem} Let $F(x,y)$ be a homogeneous polynomial of degree $n<p$ in $\mathbb{F}_{q}((x))[y]$ whose associated polynomial acts as the minimum polynomial for its roots, and let $\odot$ denote the composed product operation on $\mathbb{M}_{h,min}$.  Suppose that $F$ decomposes over $\mathbb{M}_{h,min}$ as $$F(x,y)=(f_1\odot f_2 \odot \cdots \odot f_t)(x,y)$$
 where $f_i(x,y) \in \mathbb{M}_{h,min}$ for each $i= 1, \dots, t$ and $(deg(f_{i}),deg(f_{j}))=1$ for $i \neq j$. Suppose that $F(x,y)$ can be decomposed in an alternate fashion as follows:
\begin{center}
$F(x,y)=(g_{1} \odot g_{2} \odot \cdots \odot g_{s})(x,y)$,
\end{center}
where $g_{i}(x,y) \in \mathbb{M}_{h,min}$ for each $i$ from 1 to $s$ and $(deg(g_{i}),deg(g_{j}))=1$ for $i \neq j$. Then $s=t$ and there is some reordering of the $g_{i}$'s so that $f_{i} \sim g_{i}$ for $i$ from 1 to $t$.
\end{Theorem}

\begin{proof} The Puiseux expansion of $F$ is $y=ax$ for some element $a \in \Gamma_{q}$ of degree $n$ over $\mathbb{F}_{q}$. Since $f(x,y)$ decomposes as $F(x,y)=(f_{1} \odot f_{2} \odot \cdots \odot f_{t})(x,y)$, it follows from the definition of the $\odot$-operation that $a$ can be written as $a=a_{1}a_{2} \cdots a_{t}$ where $y=a_{i}x$ is a Puiseux expansion of $f_{i}$ for each $i$. Likewise since $F(x,y)$ decomposes as $F(x,y)=(g_{1} \odot g_{2} \odot \cdots \odot g_{s})(x,y)$, it follows that $a$ can be written as $a=b_{1}b_{2} \cdots b_{s}$ where $y=b_{i}x$ is a Puiseux expansion of $g_{i}$ for each $i$. By Corollary 1.3, $s=t$ and there is a reordering of the $b_{i}$'s so that $a_{i}=c_{i}b_{i}$ for each $i$ from 1 to $t$ where $c_{i} \in \mathbb{F}_{q}$ for each $i$ and $c_{1}c_{2} \cdots c_{t}=1$. Thus for each $i$ there exists a polynomial $h_{i}(x,y) \in \mathbb{M}_{h,min,1}$ such that $f_{i}(x,y)=(h_{i} \odot g_{i})(x,y)$, that is $f_{i} \sim g_{i}$ for $i$ from 1 to $t$.
\end{proof}

\section{Summary}
We have reviewed some fundamental notions regarding both the composed
product of univariate polynomials over a finite field and Puiseux
expansions for plane curves, and have combined key aspects of each subject
to form a bivariate counterpart to the composed product operation. The
definition of the bivariate composed product is not dependent upon the
characteristic of the field from which the coefficients come. Examples of
bivariate composed products of various forms have been given, and some
properties of this new composition have been introduced, particularly for
the case in which the algebraic closure has positive characteristic. In
particular, we have shown that the decomposition theorem for the
univariate composed product has an analogue in the bivariate case for
certain homogeneous polynomials in $\mathbb{F}_{q}[x,y]$. 

The question of whether there is a suitable multivariate (here, three or
more variables) counterpart to the composed product is still open.
Gonzalez, Getino, and Farto have recently produced a multivariate
extension of the Newton-Puiseux algorithm \cite{ggf}, but the algorithm
relies on a monomial ordering compatible with the order of the terms of a
Taylor series expansion, and thus one cannot guarantee that the resulting
Puiseux expansion will be unique.  Note: The authors have partially addressed 
this matter in a subsequent paper, see \cite{MN}

\section{Acknowledgments}
The authors would like to thank Shuhong Gao (Clemson University), the referees, and the editor for their insightful comments on the manuscript.


\begin{thebibliography}{A}
\bibitem[1]{Ab}
s. Abhyankar, {\em Algebraic Geometry for Scientists and Engineers}, American Mathematical Society, Providence, 1990.

\bibitem [2] {bb} J.V. Brawley and D. Brown, {\it Composed Products and Module Polynomials over Finite Fields.\/} Discrete Mathematics {\bf 117} (1993), pp. 41-56.

\bibitem [3] {bc1} J.V. Brawley and L. Carlitz, {\it Irreducibles and the Composed Product for Polynomials Over a Finite Field.\/} 
Discrete Mathematics {\bf 65} (1987), 115-139.

\bibitem [4] {bc2} J.V. Brawley and L. Carlitz, {\it A Test for Additive Decomposability of Irreducibles Over a Finite Field.\/} Discrete Mathematics {\bf 76} (1989), pp. 61-65.

\bibitem [5] {bgm1} J.V. Brawley; S. Gao; and D. Mills, {\it Computing Composed Products of Polynomials.\/}  Finite fields: theory, applications, and algorithms (Waterloo, ON, 1997), pp. 1-15, Contemp. Math. {\bf 225}, Amer. Math. Soc., Providence, RI, 1999.

\bibitem [6] {bgm2} J.V. Brawley, S. Gao, D. Mills, {\it Associative rational functions in two variables, Finite Fields and Applications} (D. Jungnickel and H. Niederreiter, eds.), 43-56. Proceedings of the Fifth International Conference on Finite Fields and Applications. Springer-Verlag, Berlin, 2001.

\bibitem[7]{BK}
E. Brieskorn and H. Kn\"orrer, {\em Plane Algebraic Curves},
Birkh\"auser, Boston, 1986.

\bibitem[8]{Ca}
A. Campillo, {\em Algebroid Curves in Positive Characteristic}, Springer-Verlag, New York, 1980.

\bibitem[9]{ggf}
A. Gonz\'{a}lez; J. Getino,  and J. Farto, {\it An Algorithm for an Eigenvalues Problem in the Earth Rotation Theory}, Journal of Computational and Applied Mathematics, {\bf 10}, 1999, 243-254.

\bibitem [10] {mil} D. Mills, {\it Factorizations of root-based polynomial compositions}, Discrete Mathematics 240 (2001), no. 1-3, 161-173.

\bibitem[11] {MN} D. Mills, K. Neuerburg, {\it Root-based compositions of multivariate polynomials: structure, geometric interpretations, and decomposition results}, to appear in Communications in Algebra."

\bibitem[11]{Ne}I. Newton, {\em The Correspondence of Isaac Newton, volume II},
Cambridge University Press, London, 1960.

\bibitem[12]{Wa}
R. Walker, {\em Algebraic Curves}, Princeton University Press, Princeton, 1950.
\end{thebibliography}
\end{document}